\documentclass[english,a4paper]{article}

\usepackage[english]{babel}
\usepackage[ paperheight  =297mm,paperwidth   =210mm,  
             layoutheight =200mm,layoutwidth  =120mm,
             layoutvoffset= 48.9mm,layouthoffset= 45mm,
             centering,
             margin=0pt, includeheadfoot,
             footskip=8mm,
           ]{geometry}

\usepackage{booktabs,amsmath,amsfonts,amssymb}
\usepackage[mathscr]{euscript}
\usepackage{tikz}
\usetikzlibrary{arrows,decorations,shapes,automata}

\newtheorem{definicion}{Definition}

\tikzstyle{mundogrande}=[circle,draw,inner sep=0pt,minimum size=27pt]
\tikzstyle{mundo}=[circle,draw,inner sep=0pt,minimum size=20pt]
\tikzstyle{mundomedio}=[circle,draw,inner sep=0pt,minimum size=18pt]
\tikzstyle{minitexto}=[inner sep=0pt,minimum size=10pt]

\newcommand{\tupla}[1]{\langle #1 \rangle}
\newcommand{\set}[1]{\{{#1}\}}

\newcommand{\bigset}[1]{\big\{ {#1} \big\}}

\newcommand{\power}[1]{\wp(#1)}

\newcommand{\ov}[1]{\overline{#1}}
\newcommand{\bs}[1]{\boldsymbol{#1}}
\newcommand{\talque}{\bs{\mid}}

\newenvironment{ctabular}[1]
{\begin{center}\begin{tabular}{#1}}
{\end{tabular}\end{center}}

\newenvironment{ltabular}[1]
{\begin{flushleft}\begin{tabular}{#1}}
{\end{tabular}\end{flushleft}}

\newcommand{\limp}{\rightarrow}
\newcommand{\ldimp}{\leftrightarrow}

\newcommand{\mDia}[1]{\langle#1\rangle\,}

\newcommand{\pa}{\mathtt{P}}
\newcommand{\LA}{\ensuremath{\mathscr{L}}}

\newcommand{\km}{K}
\newcommand{\bm}{B}

\newcommand{\pla}{\leq}                 
\newcommand{\ind}{\sim}                 

\newcommand{\boxind}{[\ind]}
\newcommand{\diaind}{\langle\ind\rangle}
\newcommand{\mmboxind}[1]{\boxind\,{#1}}
\newcommand{\mmdiaind}[1]{\diaind\,{#1}}

\newcommand{\boxle}{[\pla]}
\newcommand{\diale}{\langle\pla\rangle}
\newcommand{\mmboxle}[1]{\boxle\,{#1}}
\newcommand{\mmdiale}[1]{\diale\,{#1}}

\newcommand{\mk}[1]{{\km}{#1}}            

\newcommand{\mb}[1]{{\bm}{#1}}            

\newcommand{\opObs}[1]{#1!}                         
\newcommand{\MOObs}[2]{{#1}_{\opObs{#2}}}           
\newcommand{\opUpg}[1]{{#1}\Uparrow}                
\newcommand{\MOUpg}[2]{{#1}_{\opUpg{#2}}}           

\begin{document}

\setcounter{page}{100}

\title{Which is the least complex explanation? Abduction and 
  complexity\thanks{Draft of the paper published in Max A. Freund, Max
    Fernandez de Castro and Marco Ruffino (eds.), \emph{Logic and
      Philosophy of Logic. Recent Trends in Latin America and
      Spain}. College Publicacions. Studies in Logic, Vol. 78, 2018,
    pp. 100--116.}} 
\author{Fernando Soler-Toscano}
\date{Grupo de L\'ogica, Lenguaje e Informaci\'on\\
  University of Seville\\
  \texttt{fsoler@us.es}}
\maketitle{}

\section{How an abductive problem arises?}
\label{sec:como-aparece-un}

What is abductive reasoning and when is it used? First, it is a kind
of inference. In a broad sense, a «logical inference» is any operation
that, starting with some information, allows us to obtain some other
information. People are continually doing inference. For example, I
cannot remember which of the two keys in my pocket opens the door of
my house. I try with the first one but it does not open. So I
conclude that it should be the other. In this case, the inference
starts with some data (premises): one of the two keys opens the door,
but the first I tried did not open. I reach some new information
(conclusion): it is the second key. 

Inference (or reasoning) does not always follow the same way. In the
example, if I know that one of the keys opens my house and I cannot
open with the first key, it is \emph{necessary} that the second
opens. When this is the case (that is, the conclusion follows
\emph{necessarily} from the premises), then we are facing a \emph{deductive
  inference}. Given that the conclusion is a necessary consequence of
the premises, there is no doubt that the conclusion is true, whenever
premises are all true: the door must be opened with the second key.  

But there are many contexts in which we cannot apply deductive
reasoning. Sometimes, we use it but we become surprised by the
outcome. What happens if finally the second key does not open the
door? This kind of \emph{surprise} was studied by the philosopher
Charles S. Peirce as the starting point of abductive reasoning: 
\begin{quote}
The surprising fact, \emph{C}, is observed;\\
But if \emph{A} were true, \emph{C} would be a matter of course,\\
Hence, there is reason to suspect that \emph{A} is true. (\emph{CP 5.189, 1903}).
\end{quote}
Peirce mentions a surprising fact, \emph{C}, that in our example is
that none of the keys opens the door, despite we strongly believed
that one of them was the right one. This is an abductive problem: a
surprising fact that we cannot explain with our current
knowledge. Then, we search for a solution, an explanation \emph{A}
that would stop \emph{C} from being surprising. To \emph{discover}
this \emph{A} we put into play our knowledge about how things usually
happen. For example, we may realise that maybe someone locked the 
door from the inside. Also, if someone usually locks the door
from the inside, then the explanation becomes stronger and, as Peirce says,
«there is reason to suspect» that it is true.

Not all abductive problems are identical. A common distinction is
between \emph{novel} and \emph{anomalous} abductive
problems~\cite{ali05}. A novel abductive problem is produced when the
surprise produced by \emph{C} is coherent with our previous
information. Contrary, in an anomalous abductive problem we previously
thought that \emph{C} could not be the case, and the surprise
contradicts our previous belief. The example of the key that does not
open the door is a case of anomalous abductive problem.

To clarify the notions, we now offer some informal definitions of the
concepts that are commonly used in the logical study of abductive
reasoning~\cite{libroABDFsoler}. Suppose that the symbol $\vdash$
represents our reasoning ability, so that $A,B \vdash C$ means that
from premises $A$ and $B$ it is possible to infer $C$ by a
\emph{necessary} inference (deduction, as explained above). The
negated symbol, as in $A,B \not\vdash C$, means that the conclusion
$C$ cannot be obtained from premises $A$ and $B$. Also, consider that
$\Theta$ is a set of sentences (logical propositions) representing our
knowledge (that I have two keys, one of them is the right one, etc.)
and $\varphi$ is the surprising fact (none of the keys opens the
door). Then, in a novel abductive problem $(\Theta,\varphi)$ the
following holds:
\begin{enumerate}
\item $\Theta\not\vdash \varphi$
\item $\Theta\not\vdash \lnot\varphi$
\end{enumerate}
The first condition is necessary for $\varphi$ to be surprising: 
it does not follow from our previous knowledge. The second condition
is specific for a novel abductive problem: the negation (the opposite)
of $\varphi$, represented by $\lnot\varphi$, does not follow 
from our knowledge $\Theta$. So, in a novel abductive problem our
previous knowledge was not useful to predict either the surprising
fact $\varphi$ or the contrary $\lnot\varphi$. 

In an anomalous abductive problem $(\Theta,\varphi)$, the conditions
that are satisfied are the following: 
\begin{enumerate}
\item $\Theta\not\vdash \varphi$
\item $\Theta\vdash \lnot\varphi$ 
\end{enumerate}
Now, although the first condition is the same, the second is
different: our previous knowledge $\Theta$ predicted $\lnot\varphi$,
the negation of the surprising fact $\varphi$.

\section{How an abductive problem is solved?}
\label{sec:como-se-resuelve}

Logicians say that in deductive reasoning the conclusion is contained
in the premises. This means that the information given by the
conclusion is implied by the information in the premises. For example,
the information that one of my keys open the door but the first does
not open contains the information that the second key will open. But
this does not happen in abductive reasoning: the information that the
door is locked from the inside is not implied by the information of my
keys not opening the door. So, abductive reasoning raises conclusions
that introduce new information not present in the premises. Because of
this, abduction requires a dose of creativity to propose the
solutions. Moreover, there are frequently several different solutions,
and the ability to select the best of them is required. We will return
later to this issue.

We have distinguished two kinds of abductive problems. Now we will
comment the kinds of abductive solutions that are usually
considered. We denoted above by $(\Theta,\varphi)$ an abductive
problem that arises when our knowledge is represented by $\Theta$ and
the surprising fact is $\varphi$. The solution to this problem is
given by some information $\alpha$ such that, together with the
previous knowledge we had, allows us to infer $\varphi$, logically
represented by 
$$\Theta,\alpha \vdash \varphi$$
This is the minimal condition for an abductive solution $\alpha$
to solve the problem $(\Theta,\varphi)$. Atocha Aliseda~\cite{ali05}
calls \emph{plain} to those abductive solutions satisfying this
requirement. 

There are other very interesting kinds of abductive solutions. For
example, \emph{consistent} solutions satisfy the additional condition
of being coherent with our previous knowledge. It is formally
represented by  
$$\Theta,\alpha\not\vdash\bot,$$
where the symbol $\bot$ represents any contradiction. It is
important that our abductive solutions are consistent. Possibly, we
will not know whether the abductive solution $\alpha$ is true, but
usually, if it is inconsistent with our previous knowledge, we have
reason to discard it.

Finally, \emph{explanatory} abductive solutions are those satisfying 
\[\alpha\not\vdash\varphi,\]
that is, the surprising fact $\varphi$ cannot be inferred with
$\alpha$ alone without using the knowledge given by 
$\Theta$. This is to avoid self-contained explanations: the key idea
behind this criterion is that a good abductive explanation offers the
missing piece to solve a certain puzzle, but all the other pieces were
previously given. 

When an abductive solution satisfies the three conditions above, we
call it a \emph{consistent explanatory solution}. To avoid useless or
trivial solutions (the key does not open the door because it does not
open it), it is frequent to focus on consistent explanatory
abduction. 

The logical study of abductive reasoning has been receiving a notable
attention for several years, and many calculi have been proposed for
abduction in different
logical systems~\cite{ciapir,libroABDFsoler,soler09jancl}. Now, we are
not interested in offering a specific calculus for a particular logic,
but in looking to an old problem in abductive reasoning: the
selection of the best hypothesis. Which is the best abductive
solution? First, we will proceed conceptually, by introducing some
notions from information theory. It will be in
Section~\ref{sec:una-propuesta-de} when, as an example, we will apply
the introduced idea in the context of epistemic logic.

\section{Which is the least complex explanation?}
\label{sec:cual-es-la}

It may happen that for a certain abductive problem there are several
possible explanations, not all of them mutually compatible. For
example, to explain why the key does not open the door, we have
proposed that someone locked it from the inside. But it could also
happen that the key or the door lock are broken, or that someone
changed the lock while we were outside, or made a joke, etc. It is
necessary to select one of the many possible explanations, because it
cannot be that case that all of them happened, it is enough just one
of them to explain that we cannot open the door with our key. What
explanation is selected and which criteria are used to select it? This
is the well-known problem of the selection of abductive
hypotheses~\cite{FerSV14jlli}.

Moreover, different to deductive reasoning, abductive conclusions
(selected solutions) are not necessary true. It is easy to observe
that, despite we think that someone locked the door from the inside, it
may have not been the case, and that in fact the lock is broken. So,
we usually have to replace an explanation with another one, when we
come to know that the originally chosen is false. 

Several criteria have been proposed to solve the problem of the
selection of abductive hypotheses. A common one is 
\emph{minimality}, that prefers explanations assuming fewer pieces of
new information. So, if I can solve a certain abductive problem both
assuming $\alpha_1$ or $\alpha_2$, it is possible that I can also
solve it by simultaneously assuming $\alpha_1$ and $\alpha_2$, or
maybe $\alpha_1$ and a certain $\beta$, but we will usually discard
those options because they are not the simplest possible ones. In
logical terms, if  $A\vdash B$ and both $A$ and $B$ can solve a
certain abductive problem, we prefer $B$, given that $A$ is at least
equally strong than $B$, and maybe stronger, in the sense of assuming
more information. 

Frequently, the minimality criterion is not enough to select the best
explanation. Which one is simpler: to think that someone locked the door from
the inside, or that they spent a joke by changing the door lock? 

Are there criteria that can help us to select the simplest explanation
in a broad spectrum of abductive problems? To give an (affirmative)
answer to this question we will move to a field in theoretical
computer science: Algorithmic Information Theory (AIT), which is due to the
works of Ray Solomonoff~\cite{Solomonoff1964p1,Solomonoff1964p2}, Andr\'ei
Kolmog\'orov~\cite{Kolmogorov65}, Leonid Levin~\cite{levin1973notion}
and Gregory Chaitin~\cite{Chaitin1975}. 

A central notion in AIT is the measure known as \emph{Kolmog\'orov
  complexity}, or \emph{algorithmic complexity}. To understand it, let
us compare these two sequences of 0s and 1s: 
\begin{center}
  0101010101010101010101010101010101010101\\
  0001101000100110111101010010111011100100
\end{center}
If we were asked which one of them is simpler, we will answer that the
first one. Why? It is built up from 20 repetitions of the pattern
«01». The second sequence is a random string. The difference between
the regularity of the first sequence and the randomness of the second
one is related with one property: the first sequence has a much
shorter description than the second. The first sequence can be
described as `twenty repetitions of «01»', while the second one can be
hardly described with a description shorter than itself. 

The idea behind the notion of Kolmog\'orov complexity is that if some
object $O$ can be fully described with $n$ bits (\emph{bit}: \emph{bi}nary
digi\emph{t}, information unit), then $O$ does not contain more
information. So the shorter description of the object $O$ indicates
how much information is contained in $O$. We would like to measure in
this way the complexity of abductive solutions, and introduce an
informational minimality criterion: we select the least complex
explanation, that is, the least informative one. But we will look at
how this complexity measure is quantified. 

Kolmog\'orov uses the concept of universal Turing
machine~\cite{turing1936a}. An universal Turing machine (UTM) $M$ is a
programmable device capable of implementing any algorithm. The
important point for us now is that the machine $M$, similar to our
computers, takes a program $p$, runs it and eventually (if the
computation stops) produces a certain output $o$. To indicate that $o$
is the output produced by UTM $M$ with program $p$ we write 
$M(p)=o$. Then, for a certain string of characters $s$, we define its
Kolmog\'orov complexity, $K_M(s)$, as
\[
   K_M(s) =\ \min\ \{ l(p) \mid M(p) = s, \ p \ \text{is a program}\}
\]
where $l(p)$ is the length in bits of the program $p$. That is,
$K_M(s)$ is equal to the size of the shorter program producing $s$ in
the UTM $M$. The subindex $M$ in $K_M(s)$ means that its value depends
on the choice of UTM, because not all of them interpret the programs
in the same way, despite all having the same computational power. If
we choose another machine $M'$ instead of $M$, it can happen that
$K_{M'}(s)$ is pretty different to $K_M(s)$. However, these bad news
are only relative, given that the Invariance Theorem guarantees that
the difference between $K_M(s)$ and $K_{M'}(s)$ is always lower than a
certain constant not depending on $s$, but on $M$ and $M'$. So, as we
face more and more complex strings, it is less relevant the choice of
UTM. Then, we can simply write $K(s)$ to denote the Kolmog\'orov
complexity of $s$. 

The use of Turing machines and programs allows us to set an encoding
to describe any computable (that is, that can be produced by some
algorithm) object $O$. Then, for the first binary sequence above, the
shortest description will not be `twenty repetitions of «01»' (26
characters) but the shortest program producing that string. 

To approach the relation between algorithmic complexity and abductive
reasoning, we can look at the work of Ray
Solomonoff, that conceives algorithmic complexity as a tool to create
a model that explains \emph{all the regularities in the observed
  universe} (see~\cite{Solomonoff1964p1}, Section 3.2). The idea of
Solomonoff is ambitious, but it is in line with the common postulates
of the inference to the best explanation~\cite{lip91}. A theory can be
conceived as a set of laws (axioms, hypotheses, etc., depending on the
kind of theory) trying to give account of a set of observations in a
given context. The laws in the theory try to explain the regularities
in those observations. So, what is the best theory? From Solomonoff's
point of view, the best theory is the most compact one, that
describing the highest number of observations (the most general one)
with the fewest number of postulates (the most elegant from a logical
point of view). It is the condition for the lowest algorithmic
complexity. The best theory is then conceived as the shortest program
\emph{generating} the observations that we want to explain. Such
\emph{generation} consist of the inferential mechanism underlying the
postulates of the theory. If it were a set of logical rules, the
\emph{execution} of the program given by the theory is equivalent to
what logicians denote by the \emph{deductive closure} of the theory:
the set of all consequences that can be deduced from the theory
axioms. Such execution does not always finish in a finite number of
steps, because logical closures frequently (always in classical logic)
are infinite sets, but usually there are procedures (in decidable
logical systems) that, in a finite number of steps, check whether a
certain formula belongs to such closure. 

As we can see, the algorithmic complexity measure $K(s)$ can be used
to determine which is the best theory within those explaining a set of
observations. However, the problem with $K(s)$ is its
uncomputability: there is no algorithm such that, given the object $s$
(a binary string or a set of observations) returns, in a finite number
of steps, the value $K(s)$ (the size of the shortest program producing
$s$, or the smallest theory explaining our observations). Therefore, if we
can not generally know the value of $K(s)$, we cannot know which
is the shortest program (or theory) generating the string $s$ (explaining
our observations). The most interesting consequence of the above
explanation is that, in general, the problem of determining which is the
\emph{best explanation} for a given set of observations is
uncomputable (if we understand \emph{the best} as \emph{the most
  compact}). 

However, despite the uncomputability of $K(s)$, there are good
approximations that allow to measure the algorithmic complexity of an
object. One of the most used approximations is based on lossless
compression algorithms. These algorithms are frequently used in our
computers to compress documents. A very common compression algorithm
is Lempel-Ziv, on which the ZIP compression format is based. It allows
to define a computable complexity measure that approximates
$K(s)$~\cite{LempelZiv76}. If we have some file $s$ and the output of
the compression algorithm is $c(s)$ (it is important to use a lossless
compression algorithm so that when decompressed it produces exactly
the original file $s$), we can understand $c(s)$ as a program that,
when is run in certain computer (the decompressor program) produces
$s$. So, the length of the compressed file $c(s)$ is an approximation
to $K(s)$. It is not necessary that $c(s)$ is the shortest possible
description of $s$, as we can consider it an approximation. In fact,
many applications based on $c(s)$ to measure complexity are used in
different disciplines line physics, cryptography or
medicine~\cite{LiVitanyi08}.

How can we approximate $K(s)$ to compare the complexity of several
abductive explanations and choose the best one? Compression-based
approximations to $K(s)$ are frequently good when $s$ is a character
string. It also happens with other approximations based on the notion
of \emph{algorithmic probability}~\cite{d5}. However, abductive
explanations are usually produced in the context of theories with a
structure that can be missed when treated as character
strings. However, we can use several tricks to reproduce some aspects of
the structure of the theories into the structure of the strings. For
example, in classical propositional logic, both sets of formulas
$A=\{p\to q\}$ and $B=\{\lnot q\to \lnot p,\ \lnot p\lor q\}$ 
are equivalent, but if we understand $A$ and $B$ as character
sequences and we compress them, $B$ will probably seem more complex
than $A$. We can avoid this problem by converting both sets into a
normal form, for example the minimal clausal form---sets
(conjunctions) of sets (disjunctions) of literals (propositional
variables or their negations)---which in both cases is $\{\{\lnot
p, q\}\}$.  

Another important point to be considered is that the complexity of an
abductive explanation $\alpha$ should be measured related to the
context in which it is proposed: the theory $\Theta$. Hence, $K(\alpha)$
may not be a good approximation to the complexity of $\alpha$ as an
abductive solution to a certain abductive problem
$(\Theta,\varphi)$. Because of that, in certain cases it is more
reasonable to use the notion of \emph{conditional algorithmic
  complexity} $K(s\mid x)$ measuring the length of the shortest
program that produces $s$ with input $x$. So, $K(\alpha\mid \Theta)$
would be a better approximation to the complexity of the abductive
solution $\alpha$ (within the theory $\Theta$) than just
$K(\alpha)$. Using lossless compression, if $c(f)$ represents
the compression of $f$ and  $|c(f)|$ is the length in bits of 
$c(f)$, a common approximation to $K(s\mid x)$ is given by
$|c(xs)|-|c(x)|$, where $xs$ represents the concatenation of $x$ and
$s$. It can be observed that, in general, this approximation gives
different values for $K(y\mid x)$ and $K(x\mid y)$, and the value of 
$K(x\mid x)$ approaches $0$ for an ideal compressor, 
given that the size of the compression of $xx$ is almost equal to the
compression of $x$, only one instruction to repeat all the output has
to be included. 

As we can see, the notions of algorithmic complexity make sense to
approach the problem of the complexity of abductive solutions and to
tackle with computational tools the problem of the selection of the
best explanation. However, good choices have to be made about the way
to represent the theories (for example, in clausal form) and which
approach to $K(s)$ o $K(s\mid x)$ is to be used.  We presented above a
very simple example on propositional logic where clausal form can be
fine. But other options are also possible. For example, Kripke frames
can be used to represent relations between
theories~\cite{Soler-Toscano01042012}. That way, each world $w$
represents a possible theory $\Theta_w$, and the accessibility
relation indicates which modifications can be done to the
theories. Then, if world $w$ can access to $u$, then $\Theta_w$ can be
modified to become $\Theta_u$. An abductive problem appears when we
are in a certain world $w$ and there is a certain formula $\varphi$
which does not follow from theory $\Theta_w$. Then, we solve the
abductive problem by moving to another accessible world $u$ (we modify
our theory $\Theta_w$ to get $\Theta_u$) such that $\varphi$ is a
consequence of $\Theta_u$. This way we give account of modifications
in theories that go beyond adding new formulas. That is, if does not
necessary happen $\Theta_w\subset \Theta_u$, because the change of
theory can entail deeper modifications, for example in the structural
properties of the logical consequence relation. Then, it may be
possible to pass, for example, from $\Theta_w$ with a monotonous
reasoning system, to a non-monotonic reasoning in $\Theta_u$. However,
within all accessible theories from $w$ that explain $\varphi$, the
problem of determining the least complex explanation still
remains. The complexity measure that should be used here is
$K(\Theta_u\mid \Theta_w)$ and, among all accessible theories from $w$
explaining $\varphi$, the one which minimises this complexity measure
should be chosen.

Despite offering resources to compare different abductive solutions
and to choose the simplest one, algorithmic
complexity notions have two problems: (1) to determine a good
representation for theories (or formulas) and (2) to choose a
computable approximation to $K(s)$ or $K(s\mid x)$, through lossless
compression or by other means. In the next section we present an
example, based on epistemic logic, illustrating how we can do this in
a specific case.

\section{A proposal using epistemic logic}
\label{sec:una-propuesta-de}

In this section we introduce an application of $K(s)$ to the selection
of the best abductive explanation, in the context of dynamic epistemic
logic (DEL). The presentation is based on previous papers where we use
the same logical tools~\cite{soler14theoria,FerSV14jlli,abdEpistIGPL},
but the selection criteria are now different. Here, an approximation
to $K(s)$ is applied to choose among several abductive explanations.

One of the possible ways to model the knowledge and belief of an agent
is offered by \emph{plausibility
  models}~\cite{BaltagSmets2008tlg}. We start by presenting the
semantic notions that will be later used to propose and solve
abductive problems. 

\begin{definicion}[Language $\LA$]\label{def:BasicFullLanguage}\index{plausibility language}
  Given a set of atomic propositions $\pa$, formulas $\varphi$ of the
  language {\LA} are given by
  \begin{center}
    \begin{tabular}{@{\qquad\qquad}c@{\;::=\;}l}
      $\varphi$ & $p \mid \lnot \varphi \mid \varphi \lor \varphi \mid \mmdiale{\varphi} \mid \mmdiaind{\varphi}$ \\
    \end{tabular}
  \end{center}
  where $p \in \pa$. Formulas of the form $\mmdiale{\varphi}$ are read as ``there is a world at least as plausible as the current one where $\varphi$ holds'', and those of the form $\mmdiaind{\varphi}$ are read as ``there is a world epistemically indistinguishable from the current one where $\varphi$ holds''. Other Boolean connectives ($\land$, $\limp$, $\ldimp$) as well as the universal modalities, $\boxle$ and $\boxind$, are defined as usual ($\mmboxle{\varphi} := \lnot \mmdiale{\lnot \varphi}$ and $\mmboxind{\varphi} := \lnot \mmdiaind{\lnot \varphi}$ for the latter).
\end{definicion}

It can be observed that the language {\LA} is like propositional logic
with two new modal connectives, $\mDia{\pla}$ and $\mDia{\ind}$,
that will allow to define the notions of belief and knowledge. These
notions will depend on a \emph{plausibility order} that the agent sets
among the worlds in the model. We now see how these models are built. 

\begin{definicion}[Plausibility model]\label{def:FullModels}\index{plausibility model}
  Let $\pa$ be a set of atomic propositions. A \emph{plausibility
    model} is a tuple $M=\tupla{W, \leq, V}$, where:
  \begin{itemize}
  \item $W$ is a non-empty set of \emph{possible worlds}
  \item ${\leq} \subseteq (W\times W)$ is a locally connected and
    conversely well-founded preorder\footnote{A relation $R \subseteq
      (W \times W)$ is \emph{locally connected} when every two
      elements that are $R$-comparable to a third are also
      $R$-comparable. It is \emph{conversely well-founded} when there
      is no infinite $\overline{R}$-ascending chain of elements in
      $W$, where $\overline{R}$, the \emph{strict} version of $R$, is
      defined as $\overline{R}wu$ iff $Rwu$ and not $Ruw$. Finally, it
      is a \emph{preorder} when it is reflexive and transitive.}, the
    agent's \emph{plausibility relation}, representing the
    plausibility order of the worlds from her point of view ($w \leq
    u$ is read as ``$u$ is at least as plausible as $w$'')
  \item $V:W \to \power{\pa}$ is an \emph{atomic valuation function},
    indicating the atoms in $\pa$ that are true at each possible
    world. 
  \end{itemize}
  A \emph{pointed plausibility model} $(M, w)$ is a plausibility model
  with a distinguished world $w \in W$.
\end{definicion}

The key idea behind plausibility models is that an agent's beliefs can
be defined as what is true \emph{in the most plausible worlds from the
  agent's perspective}, and modalities for the plausibility relation
$\leq$ will allow this definition to be formed. In order to define the
agent's knowledge, the approach is to assume that two worlds are
epistemically indistinguishable for the agent if and only if she
considers one of them at least as plausible as the other (i.e., if and
only if they are comparable via $\leq$). The \emph{epistemic
  indistinguishability relation} $\sim$ can therefore be defined as
the union of $\leq$ and its converse, that is, as $\sim \,:=\, \leq
\cup \geq$. Thus, $\sim$ is the symmetric closure of $\leq$ and hence
$\leq \;\subseteq\; \sim$. Moreover, since $\leq$ is reflexive and
transitive, $\sim$ is an \emph{equivalence} relation. This epistemic
indistinguishability relation $\sim$ should not be confused with the
\emph{equal plausibility} relation, denoted by $\simeq$, and defined
as the \emph{intersection} of $\leq$ and $\geq$, that is, $\simeq
\,:=\, \leq \cap \geq$. For further details and discussion on these
models, their requirements and their properties, the reader is
referred to \cite{BaltagSmets2008tlg,Velazquez2014delieb}. 

Now we can see how a formula is evaluated in a plausibility
model. Modalities $\diale$ and $\diaind$ are interpreted in the
standard way, using their respective relations. 

\begin{definicion}[Semantic
  interpretations]\label{def:FullSemanticInt} 
  Let $(M, w)$ be a plausibility model $M=\tupla{W, \pla, V}$ with
  distinguished world $w\in W$. By $(M, w) \Vdash
    \psi$ we indicate that the formula $\psi\in \LA$ is true in the
    world $w\in W$. Formally,
  \begin{ltabular}{@{\quad}l@{\quad{iff}\quad}l}
    $(M, w) \Vdash p$  & $p\in V(w)$, for every $p\in \pa$ \\
    $(M, w) \Vdash \lnot\varphi$  &   $(M, w) \nVdash \varphi$ \\
    $(M, w) \Vdash \varphi\land\psi$  &   $(M, w) \Vdash \varphi$ and
    $(M, w) \Vdash \psi$ \\ 
    $(M, w) \Vdash \mmdiale{\varphi}$  & there exists $u \in W$ such that
    $w \pla u$ and $(M, u) \Vdash \varphi$ \\ 
    $(M, w) \Vdash \mmdiaind{\varphi}$ & there exists $u \in W$ such that
    $w \ind u$ and $(M, u) \Vdash \varphi$ \\ 
  \end{ltabular}
\end{definicion}

In plausibility models, knowledge is defined using the
indistinguishability relation. So, an agent knows $\varphi$
in some world $w$ iff $\varphi$ is true in all worlds that cannot be
distinguished from $w$ by her, that is, all worlds considered  
epistemically possible for her. However, within those worlds there is
a plausibility order, not all of them are equally plausible for the
agent. This is relevant for the notion of belief: agent believes
$\varphi$ in a certain world $w$ iff $\varphi$ is true in \emph{the
  most plausible worlds} that are reachable from $w$. Due to the
properties of the plausibility relation, $\varphi$ is true in the most
plausible worlds iff by following the plausibility order, from some
stage we only reach
$\varphi$-worlds~\cite{BaltagSmets2008tlg}. We can express this idea
with modalities $\diale$ and $\boxle$. Formally\footnote{Operator $K$
  is commonly used for knowledge. It should not be confused with
  Kolmog\'orov complexity $K(s)$ also usually represented by $K$.},
\begin{ctabular}{ll}
  Agent \emph{knows} $\varphi$    & $\mk{\varphi}:=\mmboxind{\varphi}$
  \\ 
  Agent \emph{believes} $\varphi$ &
  $\mb{\varphi}:=\mmdiale{\mmboxle{\varphi}}$ \\ 
\end{ctabular}

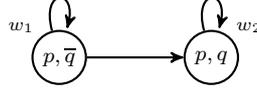
\begin{figure}[htbp!]
  \centering
    \begin{tikzpicture}[->,>=stealth',thick]
      \node [mundo] (w1)  at ( -1, 0) {\footnotesize $p, \ov{q}$};
      \node [mundo] (w2)  at (  1, 0) {\footnotesize $p, q$};

      \node [minitexto]  (t1)  at (-1.5, 0.4) {\scriptsize $w_1$};
      \node [minitexto]  (t2)  at ( 1.5, 0.4) {\scriptsize $w_2$};

      \path (w1) edge [loop above] node [midway] {} (w1)
                 edge (w2)
            (w2) edge [loop above] node [midway] {} (w2);
    \end{tikzpicture} 
  \caption{Example of a plausibility model}
  \label{fig:modExam}
\end{figure}
Fig.~\ref{fig:modExam} shows a plausibility model example
$M$. Plausibility relation $\pla$ is represented by arrows between
worlds. For the agent, world $w_2$ is more plausible than $w_1$.  In
this case, $p$ is true in both worlds, while $q$ is true only in $w_2$
($\ov{q}$ represents that $q$ is false). So, agent knows $p$ in $w_1$
but does not know $q$, that is,
$(M,w_1)\Vdash \mk{p}\land \lnot\mk{q}$. However, agent believes $q$,
$(M,w_1)\Vdash \mb{q}$. Indeed, she also believes $p$,
$(M,w_1)\Vdash\mb{p}$.
  
If a formula $\varphi$ is true in all the states of a certain model
$M$, then $\varphi$ is valid in $M$, represented as 
$M\Vdash\varphi$. In the example,  $M\Vdash \mk{p} \land \lnot
\mk{q} \land \mb{q}$. 

In the logical literature about abduction and belief revision in
general~\cite{agm85}, logical operations adding or removing
information of the theory are frequently considered. In the same way,
in the context of plausibility models, agents can perform epistemic
actions modifying the agent's information. We now present two of the
main actions that agents in plausibility models can perform. On of the
actions modifies the knowledge and the other the belief. For more
details about the properties of these actions,
see~\cite{BaltagSmets2008tlg}. 

The first operation, \emph{observation}, modifies the agent's
knowledge. It is defined in a very natural way: it consists of
removing all worlds where the observed formula is not satisfied, so
that the domain of the model is reduced.

\begin{definicion}[Observation]\label{def.update}
  Let $M=\tupla{W, \pla, V}$ be a plausibility model. The 
  \emph{observation} of $\psi$ produces the model $\MOObs{M}{\psi} =
  \tupla{W', \pla', V'}$ where
  \begin{eqnarray*}
    W' & := & \set{ w \in W \talque (M, w) \Vdash \psi }\\
    \pla' & := & {\pla} \cap {(W' \times W')}\\
    V'(w) & := & V(w)\text{, for each }  w \in W' 
  \end{eqnarray*}
\end{definicion}

This operation removes worlds of $W$, keeping only those that satisfy
(before the observation) the observed  $\psi$. The plausibility
relation is restricted to the conserved worlds. 

Another operation that agents can do is to modify just the
plausibility relation. It can be done in several ways. The operation
we call \emph{conjecture} is also known as \emph{radical upgrade} in
the literature. 

\begin{definicion}[Conjecture]\label{def:upgrade}
  Let $M=\tupla{W, \pla, V}$ be a plausibility model and $\psi$ a
  formula. The \emph{conjecture} of $\psi$ produces the model
  $\MOUpg{M}{\psi} = \tupla{W, \pla', V}$, that differs from  $M$ only
  in the plausibility relation, which is now,
  \begin{center}
    \begin{tabular}[b]{@{}l@{}c@{}l@{}}
      \;$\pla'$ & {\;:=\;} & $\bigset{ (w, u) \talque w \pla u \ \
        \text{and}\ \  (M, u) \Vdash \psi } \; \cup $ \\[1ex]
      &          & $\bigset{ (w, u) \talque w \pla u  \ \
        \text{and}\ \ (M, w) \Vdash \lnot\psi } \; \cup$ \\[1ex]
      &          & $\bigset{ (w, u) \talque w \sim u  \ \
        \text{and}\ \ (M, w) \Vdash \lnot\psi  \ \
        \text{and}\ \ (M, u) \Vdash \psi }$                      
    \end{tabular}  
  \end{center}
\end{definicion}

The new plausibility relation indicates that, after the conjecture of 
$\psi$, all $\psi$-worlds (before the conjecture) are more plausible
than all $\lnot\psi$-worlds. The previous order between 
$\psi$-worlds or between $\lnot\psi$-worlds does not
change~\cite{vanBenthem2007dlbr}. This operation preserves the
properties of the plausibility relation, as shown
in~\cite{Velazquez2014delieb}. 

We now discuss how an abductive problem can appear and be solved
within the plausibility models formalism. In the classical definition
of abductive problem, formula  $\varphi$ is an abductive problem
because it is not entailed by the theory $\Theta$. But, where does 
$\varphi$ come from? For Peirce, $\varphi$ is an observation, that is,
it comes from an agent's epistemic action. As we have seen in
Def.~\ref{def.update}, the action of observing $\varphi$ can be
modelled in DEL. What does it mean, then, that $\varphi$ is an
abductive problem? After observing $\varphi$, if it is a propositional
formula, the agent knows $\varphi$, so we cannot affirm that an
abductive problem arises when the agent does not know
$\varphi$. However, we can go back to the moment before observing
$\varphi$; if the agent did not know $\varphi$, then after the
observation it becomes an abductive problem. Formally,
\begin{equation}\label{eq:1}
  \varphi\ \ 
  \text{is an abductive problem in}\ \ 
  (\MOObs{M}{\varphi},w)\ \ 
  \text{iff}\ \ 
  (M,w)\nVdash\mk{\varphi}
\end{equation}
This definition of abductive problem within plausibility models is in
line with Peirce's idea that an abductive problem appears when the
agent observes $\varphi$. 

The notion of abductive problem in~\eqref{eq:1} has been defined in
terms of knowledge. If could have been defined in terms of belief too,
considering that $\varphi$ is an abductive problem in
$(\MOObs{M}{\varphi},w)$ iff 
$(M,w)\nVdash\mb{\varphi}$. Then, the condition for 
$\varphi$ to be an abductive problem becomes stronger than
in~\eqref{eq:1}, because $\lnot\mb{\varphi}$ implies
$\lnot\mk{\varphi}$.

Now the notion of abductive solution can also be interpreted in the
plausibility models semantics. According to Peirce's idea, the agent
knows that if $\psi$ were true, then the truth of the surprising fact
$\varphi$ would be obvious. It is expressed in DEL by requiring that
the agent knows $\psi\to\varphi$, that is,
$\mk{(\psi\to\varphi)}$. Then, when the agent faces an abductive
problem $\varphi$ and knows $\psi\to\varphi$, how does she solve it?
Again, Peirce says that 
\emph{there is reason to suspect that} $\psi$ \emph{is true}. We now
discuss what does `to suspect' $\psi$ mean in DEL and how can it be
modelled as an epistemic action. 

Something not usually considered within logical approaches to
abductive reasoning is how to integrate the solution. Maybe because in
classical logic there is no way to \emph{suspect} a formula. But in
epistemic logic there are beliefs. It is an information kind weaker
than knowledge, as we have seen. So, belief seems the most natural
candidate to model the suspicion. In line with Peirce, we then
distinguish the agent's knowledge of $\psi\to\varphi$ from her belief
in $\psi$. 

But for a reasonable suspicion, as Peirce requires, it is necessary
that the agent knows $\psi\to\varphi$. Joining all the presented
ideas, given abductive problem
$\varphi$ in $(\MOObs{M}{\varphi},w)$ (see~\eqref{eq:1}), formula
$\psi$ is a solution for it iff
\begin{equation}
  \label{eq:2}
  (M,w)\Vdash\mk{(\psi\to\varphi)}
\end{equation}
Condition
$(\MOObs{M}{\varphi},w)\Vdash\mk{(\psi\to\varphi)}$ cannot be required
because it is trivially verified in all cases in which $\varphi$ is a
propositional formula, as it becomes known after being observed, so
the agent knows also $\psi\to\varphi$ for every $\psi$. 

What does the agent do to suspect $\psi$? The most adequate abductive
action to integrate $\psi$ into the agent's information is to
\emph{conjecture} $\psi$ (def.~\ref{def:upgrade}). In this way, 
$\psi$ is integrated into the agent's information as a belief. 

\begin{figure}[htbp!]
  \centering
\begin{ctabular}{c@{}c@{}c@{}c@{}c}
  \begin{tabular}{@{\quad}c@{\quad}}
    \begin{tikzpicture}[->,>=stealth',thick]
      \node [mundo] (w1)  at ( -1, 0) {\footnotesize $p,q$};
      \node [mundo] (w2)  at (  1, 0) {\footnotesize $\ov{p},q$};
      \node [mundo] (w3)  at (  0, 1.4) {\footnotesize $\ov{p},\ov{q}$};

      \node [minitexto]  (t1)  at (-1.5, 0.4) {\scriptsize $w_1$};
      \node [minitexto]  (t2)  at ( 1.5, 0.4) {\scriptsize $w_2$};
      \node [minitexto]  (t3)  at ( -0.5, 1.8) {\scriptsize $w_3$};

      \path (w1) edge [loop above] node [midway] {} (w1)
                 edge (w2)
                 edge (w3)
            (w2) edge [loop above] node [midway] {} (w2)
                 edge (w1)
                 edge (w3)
            (w3) edge [loop above] node [midway] {} (w3)
                 edge (w1)
                 edge (w2);
    \end{tikzpicture} 
  \end{tabular}
  &
  \begin{tabular}{@{}c@{}}
    \Huge $\overset{\raisebox{4pt}{\text{\normalsize $\opObs{q}$}}}{\Longrightarrow}$
  \end{tabular}
  &
  \begin{tabular}{@{\quad}c@{\quad}}
    \begin{tikzpicture}[->,>=stealth',thick]
      \node [mundo] (w1)  at ( 0,  0.7) {\footnotesize $p,q$};
      \node [mundo] (w2)  at ( 0, -0.7) {\footnotesize $\ov{p},q$};

      \node [minitexto](t1)  at (0.3, 1.2) {\scriptsize $w_1$};
      \node [minitexto](t2)  at (0.3,-1.2) {\scriptsize $w_2$};

      \path (w2) edge [loop left] node [midway] {} (w2)
                 edge (w1)  
            (w1) edge [loop left] node [midway] {} (w1)
                 edge (w2);
    \end{tikzpicture} 
  \end{tabular}
  &
  \begin{tabular}{@{}c@{}}
    \Huge $\overset{\raisebox{4pt}{\text{\normalsize $\opUpg{p}$}}}{\Longrightarrow}$
  \end{tabular}
  &
  \begin{tabular}{@{\quad}c@{\quad}}
    \begin{tikzpicture}[->,>=stealth',thick]
      \node [mundo] (w1)  at ( 0,  0.7) {\footnotesize $p,q$};
      \node [mundo] (w2)  at ( 0, -0.7) {\footnotesize $\ov{p},q$};

      \node [minitexto](t1)  at (0.3, 1.2) {\scriptsize $w_1$};
      \node [minitexto](t2)  at (0.3,-1.2) {\scriptsize $w_2$};

      \path (w2) edge [loop left] node [midway] {} (w2)
                 edge (w1)  
            (w1) edge [loop left] node [midway] {} (w1);
    \end{tikzpicture} 
  \end{tabular}
\end{ctabular}  
  \caption{Solving an abductive problem}
  \label{fig:solAbd}
\end{figure}
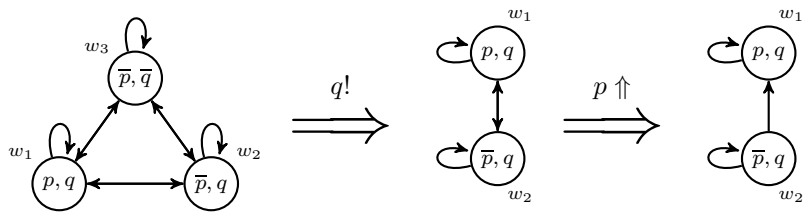
Fig.~\ref{fig:solAbd} shows an example of the whole explained process.
In the model of the left, $\mk{(p\to q)}$ is verified, but also
$\lnot\mk{q}$. In the central model, after observing $q$ agent knows
$q$, and of course she continues knowing $p\to q$, that is,
$\mk{q}\land\mk{(p\to q)}$. Then, $q$ is an abductive problem, given
that the agent knows it after the observation but not before. A
possible solution is then $p$. After the abductive action of
conjecturing $p$, the model on the right shows that the agent believes
$p$, $\mb{p}$, given that $p$ is true in the most plausible world.

Briefly,~\eqref{eq:1} shows the condition for $\varphi$ being an
abductive problem and~\eqref{eq:2} for $\psi$ being a solution for
it. As we have proposed, is it reasonable to conjecture the abductive
solution as a belief (def.~\ref{def:upgrade}).

{\medskip}

We are now ready to use Kolmog\'orov complexity as a selection criterion
within solutions of an abductive problem. First, observe that a binary
relation $R$ over a set $\{l_1,l_2,\ldots,l_n\}$ can be represented
as a binary matrix $A$ with dimension $n\times n$. In such matrix,
each cell $a_{i,j}$ is equal to $1$ iff $(l_i,l_j)\in R$ and is $0$ 
otherwise.  For example, consider relation 
$R=\{(1,1),(1,2),(2,1),(2,3),(3,1)\}$ over the set 
$\{1,2,3\}$. The matrix representing $R$ is
\[
\left(
  \begin{array}{ccc}
    1 & 1 & 0\\
    1 & 0 & 1\\
    1 & 0 & 0
  \end{array}
\right)
\]
Given a matrix $A$, its Kolmog\'orov complexity\footnote{Now $K$
  is used for algorithmic complexity and not for knowledge as in
  previous paragraphs.} $K(A)$ can be approximated for 
example using lossless compression \footnote{Lossless compression approximations are not good for
  small matrices as the one in the example, because habitual
  compressors cannot detect regularities in very short binary
  sequences. For small matrices, it is more convenient to use the
  tools presented in~\cite{d5,physA} also available at The Online
  Algorithmic Complexity Calculator:
  \texttt{http://www.complexitycalculator.com/}.}. The
least complex relations are equally the empty relation (all the matrix
is filled with 0s) and the total relation (all filled with 1s). Those
are the most compressible matrices. Other very regular relations, as
that containing only pairs $(e,e)$ for each element $e$ (diagonal
matrix), are also quite compressible. The more random a relation is,
the more random is the obtained matrix and the less compressible it is. 

In addition to approximate the value of $K(A)$ for a given matrix $A$,
compression allows us to approximate $K(B\mid A)$, given matrices
$A$ and $B$. To do that, the concatenation of $A$ and $B$ (in this
order) is first compressed and then $A$ alone is also compressed. The
difference of the size of both files approximates $K(B\mid A)$. This
approximation to $K(B\mid A)$ can be used to determine which is the
best abductive solution in the context of plausibility models. 

Consider an abductive problem $\varphi$ in 
$(\MOObs{M}{\varphi},w)$, and two competing explanations
$\psi_1$ and $\psi_2$ for it. The argument below can be extended to
any number of competing explanations. Which of both explanations will
be chosen? It was explained below that the way to integrate an
abductive explanation $\psi_i$ is to conjecture $\psi_i$ in the model
$(\MOObs{M}{\varphi},w)$ (def.~\ref{def:upgrade}). The effect of
conjecturing $\psi_i$ is only to modify the plausibility relation
of $\MOObs{M}{\varphi}$. Then, which is the best explanation? 
The answer offered by the notion of algorithmic complexity if that the
best explanation is the one modifying the model (the plausibility
relation) in the least complex way: this is not the one making the
smallest modification, but the modification with the shortest
description given the initial model. If $A_{\varphi!}$ is the matrix 
representing the plausibility relation of $\MOObs{M}{\varphi}$ and
$A_{\psi_i}$ the matrix for the plausibility relation of
$\MOUpg{\left(\MOObs{M}{\varphi}\right)}{\psi_i}$, then the best
explanation is the one minimising the value of 
\begin{equation}
  \label{eq:3}
  K(A_{\psi_i}\mid A_{\varphi!})
\end{equation}
The explanation minimising~\eqref{eq:3} is the one having the shortest
description\footnote{Formally, the chosen solution $\psi$ is the one,
  among all possible solutions, that minimises the length of a program
  that, when given to an universal Turing machine together with the
  encoding of $\MOObs{M}{\varphi}$, produces
  $\MOUpg{\left(\MOObs{M}{\varphi}\right)}{\psi_i}$.} starting at the
plausibility relation for the agent after observing $\varphi$.  

This methodology cannot be applied to small epistemic models, as those
in the previous examples, because compression is not a good
approximation to the complexity of small matrices. But it can be
applied to models of medium and large size appearing in applications
of DEL to multi-agent systems. For small models there are other
methods that can be applied~\cite{physA}, only changing in the way to
approximate $K(s)$.

\section{Discussion}
\label{sec:discusion}

Ray Solomonoff, one of the drivers of algorithmic information theory,
was convinced that the notion of algorithmic complexity can be used to
define a system explaining all the observed regularities in the
Universe. But there are strong limitations that make impossible to
build such a system, mainly the uncomputability of $K(s)$. However, by
using computable approximations, though not being able to build a
system explaining all the regularities in the Universe, it is
possible, in a specific context (plausibility models for us), to
establish criteria based on $K(s)$ (and its conditional version) that
allow to select the best explanation among the possible ones. 

There are still many issues to explore. For example, it would be
interesting to study the relevance of the notion of \emph{facticity}
introduced by Pieter Adriaans~\cite{facticity}. He considers $K(s)$
the sum of two terms, one is the structural information of $s$ and the
other the \emph{ad hoc} information in $s$. Then the best explanation
could be selected by specially looking at the contained structural
information. Also, algorithmic complexity measures can be combined
with other common selection criteria that avoid triviality. Ideally,
the least complex solution should be selected among all possible
consistent and explanatory ones.


\end{document}